\documentclass{amsart}
\usepackage[mathscr]{eucal}
\usepackage{amsfonts}
\usepackage{amsmath}
\usepackage{amsthm}
\usepackage{amssymb}
\usepackage{latexsym}
\newtheorem{theorem}{Theorem}[section]
\newtheorem{lemma}[theorem]{Lemma}
\newtheorem{proposition}[theorem]{Proposition}

\theoremstyle{definition}                               

\newtheorem{example}[theorem]{Example}

\numberwithin{equation}{section}

\def\c{\mathbb{C}}
\def\a{\mathbb{A}}
\def\k{\mathbb{K}}
\def\Z{\mathbb{Z}}

\def\N{\mathbb{N}}
\def\P{\mathbb{P}}
\def\D{\mathcal{D}}
\def\O{\mathcal{O}}
\def\h{\mathfrak{h}}
\def\g{\mathfrak{g}}
\def\gr{\mbox{\rm{Gr}}^{\mbox{\scriptsize{\rm{ad}}}}}
\def\sl2{{\mathfrak{s}\mathfrak{l}}_2}
\begin{document}
\dedicatory{To Graeme Segal for his sixtieth birthday}
\title[Differential Isomorphism]
{Differential Isomorphism and Equivalence of Algebraic Varieties}
%
%
	\author{Yuri Berest}
%
%
	\address{Department of Mathematics, Cornell University,
         Ithaca, NY 14853, USA}
	\email{berest@math.cornell.edu}
%
     \author{George Wilson}
     \address{Department of Mathematics, 
      Imperial College London SW7 2AZ, UK}
      \email{g.wilson@imperial.ac.uk}
%
%
\maketitle
\section{Introduction}
Let $\, X \,$ be an irreducible complex affine algebraic variety, 
and let $\, \D(X) \,$ be the ring of 
(global, linear, algebraic) differential operators on $ X $ 
(we shall review the definition in Section ~\ref{gen}).
This ring has a natural filtration (by order of operators)
in which the elements of order zero are just the ring $\, \O(X) \,$
of regular functions on $ \,X \,$. Thus, if we are given 
$ \, \D(X) \,$
together with its filtration, we can at once recover the variety $\, X \,$.
But now suppose we are given $\, \D(X) \,$ just as an abstract
noncommutative $\, \c $-algebra, without filtration; 
then it is not clear whether
or not we can recover $\, X \,$. We shall call two varieties  
$\, X \,$ and $\, Y \,$ 
{\it differentially isomorphic} if $\, \D(X) \,$ and
$\, \D(Y) \,$ are isomorphic. 

The first examples of nonisomorphic varieties with isomorphic 
rings of differential operators were found by Levasseur, Smith 
and Stafford (see \cite{LSS} and Section~\ref{higher} below). 
These varieties arise in the representation theory 
of simple Lie algebras;
they are still the only examples we know in dimension $\, > 1 \,$
(if we exclude products of examples in lower dimensions). 
For curves, on the other hand, 
there is now a complete classification up to differential 
isomorphism; the main purpose of this article is to review that case. 
The result is very strange. It turns out that for curves, 
$\, \D(X) \,$ determines $ \,X \,$ (up to isomorphism) 
except in the very special case when $ \,X \,$ is homeomorphic to 
the affine line $\, \a^1 $ (we call such a curve a 
{\it framed curve}).  There are uncountably many nonisomorphic framed 
curves (we can insert arbitrarily bad cusps at any finite number of 
points of $\, \a^1 $).  However, the differential isomorphism 
classes of framed curves are classified by a single non-negative 
integer $\, n \,$. This invariant $\, n \,$ seems to us the most 
interesting character in our story: it appears in many guises, some 
of which we describe in Section~\ref{n}.

We can also ask to what extent $\, X \,$ is determined by the Morita 
equivalence class of $\, \D(X) \,$: we call two varieties 
$\, X \,$ and $\, Y \,$ {\it differentially equivalent} if 
$\, \D(X) \,$ and $\, \D(Y) \,$ are Morita equivalent (as $ \c$-algebras). 
A complete classification of curves up to 
differential equivalence is not available; however, it is known that 
the differential equivalence class of a {\it smooth} affine curve 
$\, X \,$ consists of all the curves homeomorphic to $\, X \,$.  In 
particular, all framed curves are differentially equivalent to each 
other: that is one reason why the invariant $\, n \,$ which 
distinguishes them has to be somewhat unusual.  In dimension 
$\, > 1 \,$, there are already some interesting results about 
differential equivalence; we include a (very brief) survey in 
Section~\ref{higher}, where we also mention some 
generalizations of our questions to non-affine varieties.

At the risk of alienating some readers, we point out that most 
of the interest in this paper is in {\it singular} varieties.  For smooth 
varieties it is a possible conjecture that differential equivalence 
implies isomorphism: indeed, that is true for curves. 
However, in dimension $\, > 1 \,$ the conjecture would be based on no more 
than lack of counterexamples.

Our aim in this article has been to provide a readable survey, suitable as 
an introduction to the subject for beginners; most of the material is 
already available in the literature.  For the convenience of readers 
who are experts in this area, we point out a few exceptions to that rule: 
Theorem~\ref{latest} is new, and perhaps Theorem~\ref{smoothm}; also, 
the formulae (\ref{rep}) and (\ref{2}) have not previously appeared 
explicitly.

\section{Generalities on Differential Operators}
\label{gen}

We first recall the definition of (linear) differential operators, 
in a form appropriate for applications in algebraic geometry 
(see \cite{G}).  If $\, A \,$ is a (unital associative) 
commutative algebra over (say) 
$\, \c \,$, the filtered ring
$$
\D(A) = \underset{r \geq 0} \bigcup \D^{r}(A) 
\subset \mbox{\rm End}_{\c} \,(A)
$$
of differential operators on $A$ may be defined inductively as
follows. First, we set $\, \D^{0}(A) := A \,$ (here the 
elements of $\, A \,$ are identified with the corresponding  
multiplication operators); then, by definition,  
a linear map $\, \theta : A \to A \,$ belongs to 
$\, \D^{r}(A) \,$ if
$$
\theta a - a \theta \in \D^{r-1}(A) \ \, 
{\rm for\ all}\ \, a \in A \,.
$$
The elements of $ \,\D^{r}(A) \,$ are called {\it 
differential operators of order} $\, \leq r \,$ on $\, A \,$.
The commutator of two operators of orders $\, r \,$ and $\, s \,$ 
is an operator of order at most $\, r + s -1 \,$; it follows 
that the {\it associated graded algebra}
$$
\makebox{\rm gr} \, \D(A)\  := \ 
\underset{r \geq 0} \bigoplus \ 
{\D^{r}(A)}/{\D^{r-1}(A)}
$$
is commutative (we set $ \,\D^{-1}(A) := 0 \,$).

Slightly more generally, we can define the ring $\D_{A}(M)$ 
of differential operators on any $\, A $-module $\, M \,$: 
the operators of order zero are the $\, A $-linear maps $\, M \to M \,$, 
and operators of higher order are defined inductively just as 
in the special case above (where $\, M = A \,$).
\begin{example}
\label{Weyl}
If $\, A = \c[z_1 , \ldots , z_m] \,$, then 
$\, \D(A) = \c[z_i, {\partial}/{\partial z_i}] \,$ is the 
$\, m $th {\it Weyl algebra} (linear differential operators 
with polynomial coefficients).
\end{example}
\begin{example}
\label{rat}
Similarly, if $\, A = \c(z_1 , \ldots , z_m) \,$, then 
$\, \D(A) = \c(z_i) [{\partial}/{\partial z_i}]\,$ 
is the algebra of linear 
differential operators (in $\, m \,$ variables)
with rational coefficients.
\end{example}

The definition of $\, \D(A) \,$ makes sense for an arbitrary 
$\, \c$-algebra $\, A \,$; however, in this paper we shall 
use it only in the cases when $\, A \,$ is
either the coordinate ring $\, \O(X) \,$ of an 
irreducible affine variety $\, X \,$, or the 
field $\, \k \equiv \c(X) \,$ of rational functions on such a variety.
Let us consider first the latter case.  If we choose a transcendence 
basis $\,\{ z_1 , \ldots , z_m \} \,$ for $\, \k \,$ over $\, \c \,$
(where $\, m = \dim X \,$), then there are (unique) $\, \c$-derivations 
$\, \partial_1 , \ldots , \partial_m \,$ of $\, \k \,$ such that 
$\, \partial_i (z_j) = \delta_{ij} \,$, and each element of 
$ \,\D^{r}(\k) \,$ has a unique expression in the form
$$
\theta = \sum_{\vert \alpha \vert \leq r} f_{\alpha} \partial^{\alpha}
$$
(with $\, f_{\alpha} \in \k \,$), as in Example~\ref{rat} above, 
in which $\, X \,$ is the affine space $\, \a^m \,$.
In particular, $ \,\D(\k) \,$ is generated by 
$ \,\D^{1}(\k) \,$, as one would expect, and an element 
of  $ \,\D^{1}(\k) \,$ is just the sum of a derivation 
and a multiplication operator.  Indeed, it is easy to show 
that this last fact is true for an arbitrary algebra $\, A \,$. 

The case where the ring $\, A \,$ is $\, \O(X) \,$ is more subtle; 
in this case  $ \,\D(A) \,$ 
is denoted by $ \,\D(X) \,$ and is called the 
{\it ring of differential operators on $\, X \,$}. Thus the $\, m $th Weyl 
algebra (see Example~\ref{Weyl}) is the ring 
of differential operators on $\, \a^m  \,$. In general 
one does not have global coordinates on $\, X \,$, as in this example; 
nevertheless, 
if $\, X \,$ is {\it smooth}, the structure of $ \,\D(X) \,$ 
is still well understood.
\begin{proposition}
\label{smoothD}
Let $\, X \,$ be a smooth (irreducible) affine variety.  Then \\
(i) $ \,\D(X) \,$ is a simple (left and right) Noetherian 
ring without zero divisors; \\
(ii) $ \,\D(X) \,$ is generated as a $\c$-algebra by finitely 
many elements of $ \,\D^{1}(X) \,$; \\
(iii) the associated graded algebra 
$\, \makebox{\rm gr} \, \D(X) \,$
is canonically isomorphic to $\, \O(T^*X) \,$;  \\
(iv) $ \,\D(X) \,$ has global (that is, homological) 
dimension equal to $\, \dim X \,$.
\end{proposition}

If $\, X \,$ is singular, the situation is less clear.
We can still consider the 
ring $\, \Delta(X) \,$ of ($\c$-linear) operators on 
$\, \O(X) \,$ generated by the multiplication operators 
and the derivations of $\, \O(X) \,$; 
however, in general $\, \Delta(X) \,$ is 
smaller than $ \,\D(X) \,$. 
Our main reason to prefer $ \,\D(X) \,$ to 
$\, \Delta(X) \,$ is the following.  Each differential operator 
on $\, \O(X) \,$ has a unique extension to a differential operator 
(of the same order) on $\, \k \,$, so we may view $ \,\D(X) \,$ 
as a subalgebra of $ \,\D(\k) \,$.  Furthermore, a differential 
operator on $\, \k \,$ which preserves  $\, \O(X) \,$ is a  
differential operator on $\, \O(X) \,$ (this last statement 
would in general not be true for $\, \Delta(X) \,$).  Thus we have:
\begin{proposition}
\label{conc}
Let $\, X \,$ be an affine variety with function field  $\, \k \,$.  Then
$$
\D(X) = \{D \in \D(\k) : 
D.\O(X) \subseteq \O(X)\}\ .
$$
\end{proposition}

For the purposes of the present paper we could well take this as 
the definition of $ \,\D(X) \,$.  It follows from 
Proposition~\ref{conc} that $ \,\D(X) \,$ is without 
zero divisors also for (irreducible) singular varieties $\, X \,$.

\begin{example}
\label{x1}
Let $\, X \,$ be the rational curve with coordinate ring 
$\, \O(X) := \c[z^2, z^3] \,$ (thus $\, X \,$ has just 
one simple cusp at the origin).  Then $\, \Delta(X) \,$ is generated 
by $\, \O(X) \,$ and the derivations 
$\, \{ z^r \partial : r \geq 1 \} \,$ (we set 
$\, \partial := \partial / \partial z \,$).
But  
$ \,\D^{2}(X) \,$ contains the operators 
$\, \partial^2 - 2 z^{-1} \partial \,$ 
and $\, z \partial^2  - \partial \,$,
neither of which belongs to $\, \Delta(X) \,$.
\end{example}

To obtain a concrete realization of 
$\D_{A}(M)$ similar to that in Proposition~\ref{conc}, 
we need to suppose that $M$ is embedded as an 
$\, A$-submodule of some $\, \k $-vector space; to fix ideas, 
we formulate the result in the case that will concern us, where 
$\, M \,$ has rank 1.
\begin{proposition}
\label{conc2}
Suppose $\, M \subset \k\, $ is a (nonzero) 
$\, A $-submodule of $\, \k \,$. Then
$$
\D_{A}(M) = \{D \in \D(\k) : D.M \subseteq M \}\ .
$$
\end{proposition}

\subsection*{Notes}
1. To part (iii) of Proposition ~\ref{smoothD} we should add 
that the commutator on  $ \,\D(X) \,$ 
induces on $\, \makebox{\rm gr} \, \D(X) \,$ 
the canonical Poisson bracket coming from the symplectic structure 
of $\, T^*X \,$; that is, $ \,\D(X) \,$ is a 
{\it deformation quantization} of $\, \O(T^*X) \,$. 
\vspace{0.1cm} \\
2.  For singular varieties, the rings $\, \Delta(X) \,$ and 
$ \,\D(X) \,$ have quite different properties: for example, 
$\, \Delta(X) \,$ is simple if 
{\it and only if} $\, X \,$ is smooth (cf.\ Theorem~\ref{six} below).
It follows that if $\, X \,$ is smooth, then $\, \Delta(X) \,$ is 
never isomorphic, or even Morita equivalent, to $\, \Delta(Y) \,$ 
for any singular variety $\, Y \,$. Thus the present paper would 
probably be very short and dull if we were to work with $\, \Delta(X) \,$ 
rather than with $ \,\D(X) \,$. 
\vspace{0.1cm} \\
3.  Nakai (cf.\  \cite{Na}) has conjectured 
that  $ \,\D(X) = \Delta(X) \,$ if {\it and only if} 
$\, X \,$ is smooth. The conjecture has been proved for curves 
(see \cite{MV}) and, more generally, for varieties with smooth 
normalization (see \cite{T}). 
In \cite{Be} and \cite{R} it is shown that Nakai's conjecture 
would imply the well known {\it Zariski-Lipman conjecture}: if the 
module of derivations of $ \,\O(X) \,$ is projective, then 
$\, X \,$ is smooth.  
\vspace{0.1cm} \\
4.  If $\, X \,$ is singular, then in general $ \,\D(X) \,$ may have 
quite bad properties. In \cite{BGG} it is shown that if $\, X \,$ is the 
cone in $\, \a^3 \,$ with equation $\, x^3 + y^3 + z^3 = 0 \,$, then 
$ \,\D(X) \,$ is not a finitely generated algebra, nor left 
or right Noetherian. In this example $\, X \,$ is a normal variety, and 
has only one singular point (at the origin). In \cite{SS}, Section 7, 
it is shown that if $\, X \,$ is a variety of dimension $\, \geq 2 \,$ 
with smooth normalization and isolated singularities, then 
$ \,\D(X) \,$ is right Noetherian but not left Noetherian.
\vspace{0.1cm} \\
5. In the situation of Proposition~\ref{conc2}, it may happen that 
the ring $\, B := \D_{A}^{0}(M) \,$ is larger than 
$\, A \,$.  In that case the ring 
$\, \D_{A}(M) \subset \D(\k) \,$ would not change 
if we replaced $\, A \,$ by $\, B \,$; thus there is no loss of 
generality if we restrict attention to modules $\, M \,$ for which 
$\, B = A \,$.  We call such $\, A $-modules {\it maximal}.
\vspace{0.1cm} \\
6.  Of course, all the statements in this section (and, indeed, in most 
of the other sections) would remain true if we replaced $\, \c \,$ by any 
algebraically closed field of characteristic zero.  If we 
work over a field of positive characteristic, the above definition 
of differential operators is still generally accepted to be the 
correct one, but some of the properties of the rings $ \,\D(X) \,$ 
are very different: for example, $ \,\D(X) \,$ is not 
Noetherian, or finitely generated, or without zero divisors (see, 
for example, \cite{Sm}).  In particular, in positive characteristic 
$ \,\D(\a^1) \,$ is not at all like the Weyl algebra.
\vspace{0.1cm} \\
7. A convenient reference for this section is the last chapter of 
the book \cite{MR}, where one can find proofs of all the facts we have 
stated (except for Proposition~\ref{conc2}, whose proof is similar to 
that of Proposition~\ref{conc}).
\section{Differential equivalence of curves}

From now on until Section~\ref{higher}, 
$\, X \,$ will be an affine {\it curve}, probably 
singular.  In this case the problems mentioned in 
Section~\ref{gen}, Note 4 do not occur.
\begin{proposition}
\label{fg}
Let $\, X \,$ be an (irreducible) affine curve.  Then  $\, \D(X) \,$ is a 
(left and right) Noetherian ring, and is finitely generated as a 
$\, \c$-$algebra \,$.
\end{proposition}

However, the associated graded ring 
$\, \makebox{\rm gr} \, \D(X) \,$ 
is in general not a Noetherian 
ring (and hence not a finitely generated algebra either).  
The following theorem of Smith and Stafford shows that 
for our present purposes there is a very stark division of curve
singularities into ``good'' and ``bad''.
\begin{theorem}
\label{six}
Let $\, X \,$ be an affine curve, and let $\, \tilde X \,$ be its 
normalization.  Then the following are equivalent.
\begin{enumerate}
\item The normalization map $\, \pi : \tilde X \to X \,$ is bijective.
\item The algebras $\, \D(\tilde X) \,$ and $\, \D(X) \,$ 
are Morita equivalent.
\item The ring $\, \D(X) \,$ has global dimension $\, 1 \,$ 
(that is, the same as $\, \D(\tilde X) \,$).
\item The ring $\, \D(X) \,$ is simple.
\item The algebra $\, \makebox{\rm gr} \,\D(X) \,$ is 
finitely generated.
\item The ring $\, \makebox{\rm gr} \,\D(X) \,$ is 
Noetherian.

\end{enumerate}
\end{theorem}

Perhaps the most striking thing about 
Theorem~\ref{six} is that the ``good'' singularities (from our present point 
of view) are the {\it cusps} (as opposed to double points, or higher 
order multiple points). If $\, X \,$ has even one double point, the ring 
$\, \D(X) \,$ is somewhat wild; whereas if $\, X \,$ has only cusp 
singularities, no matter how ``bad'', then $\, \D(X) \,$
is barely distinguishable from the ring of differential operators on the 
smooth curve $\, \tilde X \,$.

Theorem~\ref{six} does not address the question of when two {\it smooth} 
affine curves are differentially equivalent.  However, the answer to that 
is very simple.
\begin{theorem}
\label{smoothm}
Let  $\, X \,$ and $\, Y \,$ be smooth affine curves.  Then 
$\, \D(X) \,$ and $\, \D(Y) \,$ are Morita 
equivalent (if and) only if $\, X \,$ and $\, Y $ are isomorphic.
\end{theorem}

Theorems~\ref{six} and \ref{smoothm} together determine completely 
the differential equivalence class of a smooth curve $\, X \,$: it 
consists of all curves obtained from $\, X \,$ by pinching a finite 
number of points to (arbitrarily bad) cusps.  
\subsection*{Notes}
1. Apparently, not much is known about the differential equivalence 
class of a curve with multiple points. 
From Theorem~\ref{six} one might guess that 
if  $\, \pi : Y \to X \,$ is regular surjective of degree one, then 
$\, X \,$ and $\, Y \,$ are differentially equivalent if and 
only if $\, \pi \,$ is bijective.  However, in \cite{SS} (5.8) there is  
a counterexample to the ``if'' part of this statement. 
The paper~\cite{CH2} contains some curious results about the Morita 
equivalence class of  $\, \D(A) \,$ when $\, A \,$ is the 
{\it local} ring at a multiple point of a curve. 
\vspace{0.1cm} \\
2.  Another natural question that is not 
addressed by Theorem~\ref{six} is: what is the global dimension of
$\, \D(X) \,$ if $\, X \,$ has multiple points?  In 
\cite{SS} it is proved that if the singularites are all 
{\it ordinary} multiple points, then the answer is 2; but 
for more complicated singularities it seems nothing is known. 
\vspace{0.1cm} \\
3.  We have not found Theorem~\ref{smoothm} stated explicitly 
in the literature, but it is an easy consequence of the results 
of \cite{CH1} and \cite{M-L}: we will sketch a proof in 
Section~\ref{corr}, Note 5.
\vspace{0.1cm} \\
4. Proposition~\ref{fg} is proved in \cite{SS} and (also in the case 
of a reducible (but reduced) curve) in  \cite{Muh}.
\vspace{0.1cm} \\
5.  We refer to \cite{SS} for the proofs of the  
various assertions in Theorem~\ref{six}. Here we mention only that a 
key role is played by the space
\begin{equation}
\label{P}
P \equiv \D(\tilde X , X) := 
\{ D \in \D(\k) : D.{\O(\tilde X)}
\subseteq \O(X) \} \, .
\end{equation}
Clearly, $\, P \,$ is a right ideal in $\, \D(\tilde X) \,$ 
and a left ideal in $\, \D(X) \,$; the Morita equivalence 
in Theorem~\ref{six} is defined by tensoring with the bimodule 
$\, P \,$.  Another 
notable property of $\, P \,$ is the following: each of the 
statements in Theorem~\ref{six} is equivalent to the condition
\begin{equation}
\label{P1}
P.\O(\tilde X) = \O(X) \, .
\end{equation}
The formulae (\ref{P}) and (\ref{P1}) provide the starting 
point for the theory of Cannings and Holland 
which we explain in Section~\ref{corr}\,; there  $\, P \,$ is 
replaced by an arbitrary right ideal in $\, \D(\tilde X) \,$.

\section{Differential Isomorphism of Curves}

We now turn to our main question, concerning differential 
isomorphism. We begin by sketching the history 
of this subject.

To our knowledge, the papers \cite{St}, \cite{Sm} are the first 
that explicitly pose the question: does   
$ \, \D(X) \simeq \D(Y )\,$ imply $\, X \simeq Y \,$? 
In \cite{St}, Stafford proved that this is true 
if $\, X \,$ is the affine line 
$\, \a^1 \,$ (in which case $ \,\D(X)\, $ 
is the Weyl algebra), and also if $\, X \,$ is the plane curve 
with equation $\, y^2 = x^3 \,$, that is, the rational curve 
obtained from  $\, \a^1 \,$ by introducing a 
{\it simple cusp} at the origin.  The first general result in 
the subject is due to L. Makar-Limanov (see \cite{M-L}).  His idea 
was as follows.  Recall that if we take the commutator
$\, (\mbox{ad} f) L := f L - L f \,$ 
of a function $\, f \in \O(X)\, $
with an operator $\, L \in \D(X)\, $ of order $ n $, then 
we get an operator of order at most $ n-1 $ (indeed, this is essentially 
the definition of $\, \D(X)\, $, see Section~\ref{gen} 
above).  It follows that $\, (\mbox{ad} f)^{n+1} L = 0  \,$, 
so that $\,f\,$ is a (locally) {\it ad-nilpotent} element of 
$\, \D(X)\, $.  If it happens (as seems likely) that 
the set $\, {\mathcal N}(X) \,$ of all ad-nilpotent elements of 
$\, \D(X) \,$ coincides with 
$\, \O(X) \,$, then we have a purely ring-theoretical
description of $\, \O(X) \subset \D(X) \,$, 
namely, it is the unique maximal abelian ad-nilpotent 
subalgebra (for short: {\it mad subalgebra}) 
of $\, \D(X) \,$.  So in this way  
$\, \D(X) \,$ determines $\, X \,$.   
Makar-Limanov's main remark was the following.
\begin{lemma}
\label{bisp}
Let $\, \k \,$ be the function field of a curve, and let 
$\, D \in \D(\k) \,$ have positive order. 
Let $\, N \subset \k \,$ be the 
set of elements of $\, \k \,$ on which $\, D \,$ acts ad-nilpotently. 
Then there is an element $\, q \,$ in some finite extension field 
of  $\, \k \,$ such that $\, N \subseteq \c[q] \,$.
\end{lemma}

If now $\, X \,$ is a curve such that 
$\, {\mathcal N}(X) \not= \O(X) \,$, 
that is, such that $\, {\mathcal N}(X) \,$ contains an operator of 
positive order, then it follows from Lemma~\ref{bisp} that 
$\, \O(X) \subseteq \c[q] \,$ for suitable $\, q \,$. 
Equivalently:
\begin{theorem}
\label{mak}

If $\, {\mathcal N}(X) \not= \O(X) \,$, then the 
normalization $\, \tilde X \,$ of $\, X \,$ is isomorphic to
$\, \a^1 \,$.
\end{theorem}

In his thesis (see \cite{P1}), P.~Perkins refined this result.
\begin{theorem}
\label{perk}
Let $\, X \,$ be an affine curve.  Then 
$\, {\mathcal N}(X) \not= \O(X) \,$
if and only if 
\begin{verse}
$ (i) $\ $ \tilde X \,$ is isomorphic to $ \a^1 $\,; and 

$ (ii) $\ the normalization map $\, \pi: \tilde X \to X \,$ is bijective.
\end{verse}
\end{theorem}

In other words, the differential isomorphism 
class of a curve $\, X \,$ consists just of (the 
class of) $\, X \,$ itself, except, possibly, when $\, X \,$ 
has the properties $ (i) $ and $ (ii) $ above.

For short, we shall call a curve with these two properties 
a {\it framed curve}.  More precisely, by a framed curve we shall mean 
a curve $\, X \,$ together with  a regular bijective map
$\, \pi: \a^1 \to X \,$: the choice of ``framing'' (that is,
of the isomorphism $\, \tilde X \simeq \a^1 \,$) 
is fairly harmless,
because any two choices differ only by an automorphism
$\, z \mapsto a z + b \,$ of \,$ \a^1 $. 
The two curves considered by Stafford are certainly framed curves:
Stafford's results do not contradict those of Perkins, because 
although the rings $\, \D(X) \,$ in these examples have many 
ad-nilpotent elements not in $\, \O(X) \,$, their
mad subalgebras are all 
isomorphic, so we can still extract $\, \O(X) \,$ 
(up to isomorphism) from $\, \D(X) \,$.  
For a while it might have seemed likely that the situation 
is similar for any framed curve; but counterexamples were 
found by Letzter \cite{L} and by Perkins \cite{P2}.
The following example of Letzter is perhaps the simplest 
and most striking.  Let $\, X \,$ and $\, Y \,$ be the curves 
with coordinate rings 
$$
\O(X) = \c + z^4 \c[z] \ ; \quad
\O(Y) = \c[z^2, z^5] \ .
$$
Each of $\, X \,$ and $\, Y \,$ is obtained from $\, \a^1 \,$ 
by introducing a single cusp at the origin; $\, X \,$ and $\, Y \,$ 
are clearly not isomorphic.  Indeed, we have 
$\, \O(X) \subset \O(Y) \,$, so the singularity 
of $\, X \,$ is strictly ``worse'' than that of $\, Y \,$.  
Nevertheless, Letzter proved that  $\, X \,$ and $\, Y \,$ are 
differentially isomorphic.  This example, and others in \cite{P2}, 
\cite{L}, shows that the problem of classifying {\it framed} curves 
up to differential isomorphism is nontrivial.

This problem was solved completely in the thesis \cite{K} of 
K. Kouakou.
The simplest way to state his result is as follows.
For each $\, n \geq 0 \,$, let $\, X_n \,$ denote the curve with coordinate 
ring 
\begin{equation}
\label{Xn}
\O(X_n) := \c + z^{n+1} \c[z] \ .
\end{equation}
(Thus the curves considered by Stafford are 
$\, X_0 \equiv \a^1  \,$ and $\, X_1 \,$, while the curve 
$\, X \,$ in Letzter's example above is $\, X_3 \,$).
\begin{theorem}[Kouakou]
\label{T1}
Every framed curve $\, X \,$ is differentially 
isomorphic to one of the above curves $\, X_n \,$.
\end{theorem}

On the other hand, Letzter and Makar-Limanov (see \cite{LM}) 
have proved the following.
\begin{theorem}
\label{lml}
No two of the curves $\, X_n \,$ are differentially
isomorphic to each other.
\end{theorem}

It follows that each framed curve $\, X \,$ is differentially 
isomorphic to exactly one of the special curves $\, X_n \,$: we 
shall call this number $\, n \,$ the 
{\it differential genus} of $\, X \,$, and denote it by $\, d(X) \,$. 

\subsection*{Notes}
1.  Of course, this is very unsatisfactory as a {\it definition} 
of the differential genus, because it does not make sense until after we have 
proved the two nontrivial Theorems \ref{T1} and \ref{lml}.  In 
Section~\ref{n} we discuss several more illuminating ways to 
define $\, d(X) \,$. We use the term ``genus'' because $\, d(X) \,$ 
is in some ways reminiscent of the arithmetic genus 
of a curve: it turns out that it is a sum of local contributions 
from each singular point, so it simply counts the cusps of our 
framed curve with appropriate weights.  In Section ~\ref{n} 
we shall explain how to calculate these weights: 
here we just mention that the weight of a  
{\it simple} (that is, of type $\, y^2 = x^3 \,$)
cusp is equal to 1, so if all the cusps of $\, X \,$ are 
simple, then $\, d(X)\,$ is just the number of cusps.
\vspace{0.1cm} \\
2.  Recall from Theorem~\ref{six} that the algebras $\, \D(X) \,$
(for $\, X \,$ a framed curve) are all Morita equivalent to each 
other: thus the invariant $\, d(X) \,$ that distinguishes them 
must be fairly subtle. 
\vspace{0.1cm} \\
3. Makar-Limanov's Lemma~\ref{bisp} (in a slightly disguised 
form) plays a basic role also in the theory of bispectral 
differential equations (compare the proof in \cite{M-L} with similar 
arguments in \cite{DG} or \cite{W1}).
\vspace{0.1cm} \\
4. There is no convenient reference where the reader can find a 
complete proof of Kouakou's theorem: Kouakou's thesis has never 
been published, and the (different) proof in \cite{BW1} is mostly 
omitted.  The proof that we shall explain in the next three 
sections amplifies the sketch given in \cite{W3}: it
is not the most elementary possible, but it seems to us the 
most natural available at present.  
\section{The Adelic Grassmannian}
\label{grad}

It is actually easier to prove a more general theorem than 
Theorem~\ref{T1}, as follows.  
Let $\, X \,$ be a framed curve, and let $\, \mathcal{L} \,$ 
be any rank 1 torsion-free coherent sheaf over $\, X \,$;
it corresponds to a rank 1 torsion-free $\, \O(X) $-module 
$\, M \,$.  
Then we have the ring 
$\, \D_{\mathcal L}(X) \equiv 
\D_{\O(X)}(M) \,$ 
of differential operators on (global) sections of 
$\, \mathcal{L} \,$.  If  $\, \mathcal{L} = \O_{X} \,$ 
is the sheaf of regular functions on $\, X \,$, then 
$\, \D_{\mathcal L}(X) \,$ is just the ring 
$\, \D(X) \,$ discussed previously.  Generalizing 
Theorem~\ref{T1}, we have the following.
\begin{theorem}
\label{T2}
Every algebra $\, \D_{\mathcal L}(X) \,$ is isomorphic 
to one of the algebras $\, \D(X_n) \,$.

\end{theorem}

Of course, Theorem~\ref{lml} shows that the integer $\, n \,$ in this 
assertion is unique: we call it the 
{\it differential genus of the pair}
$\, (X , \mathcal{L}) \,$ and denote it by $\, d_{\mathcal{L}}(X) \,$.

The reason Theorem~\ref{T2} is easier to prove than Theorem~\ref{T1} 
is that the space of pairs $\, (X , \mathcal{L}) \,$ has a 
large group of symmetries that preserves the isomorphism class 
of the algebra $\, \D_{\mathcal L}(X) \,$ (but does not 
preserve the subset of pairs of the form 
$\, (X , \O_{X}) \,$).  In fact the isomorphism 
classes of these pairs form the 
{\it adelic Grassmannian} $\, \gr $, a well-studied space that 
occurs in at least two other contexts, 
namely, in the theory of 
the Kadomtsev-Petviashvili hierarchy (cf.\ \cite{Kr}) and in the problem 
of classifying bispectral differential operators (see \cite{DG}, \cite{W1}).
The adelic Grassmannian is a subspace\footnote{It is perhaps the 
most interesting Grassmannian not mentioned explicitly in \cite{SW}\,.} 
of the much larger Grassmannian Gr studied in \cite{SW}.
We recall the definition of $\, \gr \,$.  For each $\, \lambda \in \c \,$, 
we choose a {\it $\lambda$-primary} subspace of 
$\, \c[z] \, $, that is, a linear subspace  $\, V_{\lambda} \,$
such that 
$$
(z - \lambda)^{N} \c[z] \subseteq V_{\lambda} \ \text{for some} \ N \, .
$$
We suppose that $\, V_{\lambda} = \c[z] \, $ for all but finitely 
many $\, \lambda \,$.  Let $\, V = \bigcap_{\lambda} V_{\lambda} \,$ 
(such a space $\, V \,$ is called {\it primary decomposable}) and, 
finally, let
$$
W = \prod_{\lambda} (z - \lambda)^{- k_{\lambda}} \, V \subset 
\c(z) \, ,
$$
where $\, k_{\lambda} \,$ is the codimension of $\, V_{\lambda} \,$ 
in $\, \c[z] \,$. By definition, $\, \gr \,$ consists of all 
$\, W \subset \c(z) \,$ obtained in this way.  The correspondence 
between points of $\, \gr \,$ and pairs $\, (X , \mathcal{L}) \,$ 
is a special case of the construction explained in \cite{SW}.
Given $\, W $, we obtain $\, (X , \mathcal{L}) \,$ by setting
$$
\O(X) := \{ f \in \c[z] : fW \subseteq W \} \, ;
$$
and $\, W \,$ is then the rank 1 $\, \O(X)$-module 
corresponding to $\, \mathcal{L} \,$.  Conversely, given 
$\, (X , \mathcal{L}) \,$, we let $\, W \,$ be the space of global 
sections of $\, \mathcal{L} \,$, regarded as a subspace of $\, \c(z) \,$ 
by means of a certain distinguished rational trivialization of 
$\, \mathcal{L} \,$ (implicitly described above).
\begin{proposition}
\label{wxl}
This construction defines a bijection between $\, \gr $ and the set 
of isomorphism classes of pairs $\, (X , \mathcal{L}) \,$, where 
$\, X \,$ is a framed curve and $\, \mathcal{L} \,$ is a maximal 
rank 1 torsion-free sheaf over $\, X \,$.
\end{proposition}

``Maximal'' here means that the $\, \O(X) $-module corresponding 
to $\, \mathcal{L} \,$ is maximal in the sense of Note 5, Section~\ref{gen}. 
\begin{example}
\label{wn}
If $\, X_n \,$ is the curve defined by (\ref{Xn}), then 
$\, \O(X_n) \,$ is $ 0 $-primary, and the corresponding 
point of $\, \gr \,$ is $\, W_n = z^{-n}\O(X_n) \,$. 
More generally, let $\, \Lambda \subset \N \,$ be any (additive) 
semigroup obtained from $\, \N \,$ by deleting a finite number 
of positive integers, and let   $\, \O(X) \,$ be the subring of 
$\, \c[z] \,$ spanned by $\, \{ z^i : i \in \Lambda \} \,$.  
Such a curve $\,  X \,$ is called 
a {\it monomial curve}; the corresponding 
point of $\, \gr \,$ is $\, z^{-m} \O(X) \,$, 
where $\, m \,$ is the number of elements of 
$\, \N \setminus \Lambda \,$.
\end{example}
\begin{example}
If $\, X \,$ has simple cusps at the (distinct) points 
$\, \lambda_1, \ldots , \lambda_r \in \c \,$, then $\, \O(X) \,$ 
consists of all polynomials whose first derivatives vanish at these 
points, and the corresponding point of $\, \gr \,$ is 
$$
W = \prod_{i=1}^{r} (z - \lambda_i)^{-1} \, \O(X) \, .
$$
More generally, if in addition 
we choose  $\, \alpha_1, \ldots , \alpha_r \in \c \,$, then 
$$
V = \{ f \in \c[z] : f'(\lambda_i) = \alpha_i f(\lambda_i) 
\  \text{for} \  1 \leq i \leq r \}
$$
is primary decomposable, and the corresponding point of $\, \gr \,$ is 
$$
W = \prod_{i=1}^{r} (z - \lambda_i)^{-1} \, V \, .
$$
In the pairs $\, (X, \mathcal{L}) \,$ here, the 
curve $\, X \,$ is the same as before, and 
as we vary the parameters $\, \alpha_i \,$ we get the various 
line bundles $\, \mathcal{L} \,$ over $\, X \,$.
\end{example}

The rings $\, \D_{\mathcal{L}}(X) \,$ that interest us are 
easy to describe in terms of $\, \gr \,$. 
If $\, W \in \gr \,$, we define the {\it ring of differential 
operators on $\, W \,$} by
$$
\D(W) := \{ D \in \c(z)[\partial] : D.W \subseteq W \}
$$
(as in Section~\ref{gen}, the dot denotes the natural 
action of differential operators on functions).
Proposition~\ref{conc2} shows:
\begin{proposition}
\label{dw}
Let $\, W \in \gr \,$ correspond to the pair $\, (X , \mathcal{L}) \,$ 
as in Proposition~\ref{wxl}.  Then there is a natural identification
$$
\D(W) \simeq \D_{\mathcal{L}}(X) \, .
$$
\end{proposition}

It remains to discuss the symmetries of  $\, \gr \,$. Some of them
are fairly obvious. 
First, we have the commutative group $\, \Gamma \,$ of the 
{\it KP flows}: it corresponds to the action 
$\, (X, \mathcal{L}) \mapsto (X, L \otimes \mathcal{L}) \,$ 
of the Jacobian (that is, the group of line bundles $\, L \,$ 
over $\, X \,$) on the space of pairs $\, (X, \mathcal{L}) \,$.
If $\, W^{\mbox{\scriptsize{\rm{an}}}} \supset W \, $ is the 
space of analytic sections 
of $\, \mathcal{L} \,$, then  $\, \Gamma \,$ is the group of maps of 
the form $\, W^{\mbox{\scriptsize{\rm{an}}}} 
\mapsto e^{p(z)} W^{\mbox{\scriptsize{\rm{an}}}} \,$, where $\, p \,$ 
is a polynomial.  Another fairly evident symmetry is the 
{\it adjoint involution} $\, c \,$ defined by
$$
c(W) = \{ f \in \c(z) : \mbox{\rm res}_{\infty} f(z) g(z) dz = 0 
\ \text{for all} \ g \in W \} \, .
$$
Like the KP flows, $\, c \,$ is just the restriction to $\, \gr \,$ of a 
symmetry of the Grassmannian Gr of \cite{SW}.  A more elusive 
symmetry of $\, \gr $ is 
the {\it bispectral involution}  $\, b \,$ 
introduced in \cite{W1}; it does not make sense on Gr\,, and does 
not have a simple description in terms of the pairs 
$\, (X, \mathcal{L}) \,$.  It can be  characterized by the formula
$$
\psi_{bW}(x,z) = \psi_{W}(z,x) \, ,
$$
where $\, \psi \,$ is the {\it stationary Baker function} 
of $\, W \,$ (see, for example, \cite{SW}). 
Let $\, \varphi = bc \,$, and let 
$\, G \,$ be the group of symmetries of $\, \gr \,$
generated by $\, \Gamma \,$ and $ \varphi \, $.  
In view of Proposition~\ref{dw},
Theorems~\ref{T1}, \ref{lml} and \ref{T2}
are all consequences of
\begin{theorem}
\label{mainth}
(i)  Let $\, V, W \in \gr \,$.  Then 
$\, \D(V) \,$ and 
$\, \D(W) \,$ are isomorphic if and only 
if $\, V \,$ and $\, W \,$ 
belong to the same $\, G $-orbit in $\, \gr \,$. \\ 
(ii)  Each orbit contains 
exactly one of the points $\, W_n \,$ from Example~\ref{wn}.
\end{theorem}

Although it is possible to formulate a proof of Theorem~\ref{mainth} 
within our present context, the proof will appear more natural if we 
use two alternative descriptions of $\, \gr \,$: we explain these 
in the next sections. First, in Section~\ref{corr} we shall see that 
$\, \gr \,$ can be identified with the space of ideals in the Weyl 
algebra $\, \D(\a^1) \,$: the ring $\, \D(W) \,$ 
then becomes the endomorphism ring of the corresponding ideal, 
and $\, G \,$ becomes the automorphism group of the Weyl algebra.  
Part (i) of Theorem~\ref{mainth} then turns into a theorem of Stafford 
(see \cite{St}).  In Section~\ref{cm} we explain how $\, \gr \,$ 
decomposes into the union of certain finite-dimensional varieties 
$\, \mathcal{C}_n \,$ that have a simple explicit description in 
terms of matrices; part (ii) of Theorem~\ref{mainth} then follows from the 
more precise assertion that these spaces  $\, \mathcal{C}_n \,$ are 
exactly the $\, G $-orbits.  
Since the action of $\, G \,$ also has a 
simple description in terms of matrices, part (ii) of the Theorem 
becomes a problem in linear algebra.
\subsection*{Notes}
1.  The fact that the action of $\, \Gamma \subset G \,$ preserves 
the isomorphism class of $\, \D(W) \,$ is almost 
trivial.  Indeed, if $\, g \in \Gamma \,$ is given (as above) by 
multiplication by $\, e^{p(z)} \,$, then  
$\, \D(gW) = e^{p(z)} \D(W) e^{-p(z)} \,$.  
It follows that $\, \D(gW) \,$ is even isomorphic to 
$\, \D(W) \,$ as a {\it filtered} algebra.  Thus 
the (filtered) isomorphism class of $\, \D_{\mathcal L}(X) \,$ 
depends only on the orbit of the Jacobian of $\, X \,$ in the space of 
rank 1 torsion-free sheaves; for example, if $\, \mathcal{L} \,$ is 
locally free, then $\, \D_{\mathcal L}(X) \,$ is isomorphic 
to $\, \D(X) \,$.  A direct proof that $\, \varphi \,$ 
preserves the isomorphism class of $\, \D(W) \,$ is also not too 
difficult: it follows from the facts that $\, \D(bW) \,$ 
and $\, \D(cW) \,$ are {\it anti}-isomorphic to 
$\, \D(W) \,$ (cf.\ \cite{BW2}, Sections 7 and 8). 
We regard the main assertions in Theorem~\ref{mainth} to be part (ii) 
and the ``only if'' statement in part (i).
\vspace{0.1cm} \\
2.  The spaces $\, W_n \,$ are fixed by $\, b \,$, so  $\, b \,$ 
induces an involutory anti-automorphism on each of the rings 
$\, \D(X_n) \,$.  Thus Theorem~\ref{T2} shows that 
the distinction between isomorphism and anti-isomorphism in 
the preceding note was immaterial.
\vspace{0.1cm} \\
3.  If  $\, \mathcal{L} \,$ is not  locally free, 
then in general $\, \D_{\mathcal L}(X) \,$ is not 
isomorphic to $\, \D(X) \,$ (see Example~\ref{34} below).
\vspace{0.1cm} \\
4.  Details of the proof of Proposition~\ref{wxl} can be found in \cite{W1}; 
see also \cite{CH4}, 1.4 and \cite{Elb}, p.\ 945.

\section{The Cannings-Holland correspondence}
\label{corr}

In this section we explain a different realization of $\, \gr \,$ 
(due to Cannings and Holland) as the space of ideals in the Weyl algebra. 
Let $\, A := \c[z, \partial] \,$ from now on denote the (first) Weyl algebra, 
and let $\, \mathcal{I} \,$ be the set of nonzero right ideals of 
$\, A \,$. Let $\, \mathcal{S} \,$ be the set of all linear subspaces 
of $\, \c[z] \,$.  If $\, V, W \in \mathcal{S} \,$, (or, later, also
if $\, V \,$ and $\, W \,$ are subspaces of $\, \c(z) \,$) we set 
\begin{equation}
\D(V,W) := \{ D \in \c(z)[\partial] : D.V \subseteq W \} \, .
\end{equation}
We define maps
$\, \alpha : \mathcal{S} \to \mathcal{I} \,$ and 
$\, \gamma : \mathcal{I} \to \mathcal{S} \,$ 
as follows. If $\, V \in \mathcal{S} \,$, we set
\begin{equation}
\label{alpha}
\alpha(V) := \D(\c[z],V) \,;
\end{equation}
and if $\, I \in \mathcal{I} \,$, we set
\begin{equation}
\label{gamma}
\gamma(I) := \{ D.\c[z] : D \in I \} \, .
\end{equation}
\begin{theorem}
\label{ch}
(i) We have $\, \alpha \gamma (I) = I \,$ if and only if 
$ \, I \cap \c[z] \not= \{0\} \,$. \\
(ii) We have $\, \gamma \alpha (V) = V \,$ if and only if
$\, V \,$ is primary decomposable. \\
(iii) The maps $\, \alpha \,$ and $\, \gamma \,$ define inverse bijections 
between the set of primary decomposable subspaces of $\, \c[z] $ 
and the set of right ideals of $\, A \,$ that intersect 
$\, \c[z] \,$ nontrivially. \\
(iv) If $\, V \,$ and $\, W \,$ are primary decomposable and 
$\, I := \alpha(V) \,$, $\, J := \alpha(W) \,$ are the corresponding 
(fractional) ideals, then 
$$
\D(V,W) = \{ D \in \c(z)[\partial] : DI 
\subseteq J \}  \simeq \mbox{\rm Hom}_{A}(I, J) \, .
$$
\end{theorem}
\begin{example}
Let $\, I_n \,$ be the right ideal
$$
I_n \ := \ z^{n+1} A \ +  \ \prod_{r=1}^{n} (z \partial - r) \, A \, .
$$
The second generator kills $\, z, z^2, \ldots, z^n \,$, so we find that 
$\, \gamma(I_n) = \O(X_n) \,$.
\end{example}

The assertions (iii) and (iv) in Theorem~\ref{ch} 
follow at once from (i) and (ii).  Now, not 
every right ideal of $\, A \,$ intersects $\, \c[z] \,$ nontrivially; 
but every ideal is isomorphic (as right $\, A $-module) to one with 
this property (see \cite{St}, Lemma 4.2).  
Furthermore, two such ideals $\, I, J \,$ 
are isomorphic if and only if $\, pI = qJ \, $ for some polynomials 
$\, p(z), q(z) \,$.  On the other hand, two primary decomposable subspaces 
$\, V, W \,$ determine the same point of $\, \gr \,$ if and only if 
$\, pV = qW \, $ for some polynomials $\, p(z), q(z) \,$; and the bijections   
$\, \alpha \,$ and $\, \gamma \,$ are clearly compatible with 
multiplication by polynomials.  
Let $\, \mathcal{R} \,$ denote the set of isomorphism classes of nonzero 
right ideals of $\, A \,$ (equivalently, of 
finitely generated torsion-free rank 1 right $\, A $-modules). Combining 
the remarks above with Theorem~\ref{ch}, we get the following.
\begin{theorem}
\label{chgr}
(i) The maps defined by the formulae (\ref{alpha}) and (\ref{gamma}) 
define inverse bijections
$$
\alpha : \gr \to \mathcal{R} \ \ \text{\rm and} \ \ 
\gamma :  \mathcal{R} \to \gr \,.  
$$
(ii)  For $\, V, \, W \in \gr \,$, there is a natural identification
$$
\D(V,W) \simeq \mbox{ \rm Hom}_A(\alpha(V), \alpha(W)) \,. 
$$
\end{theorem}

As a special case of (ii), we see that if $\, W \in \gr \,$ and 
$\, I := \alpha(W) \,$ is the corresponding ideal in $\, A \,$, 
then the algebra $\,\D(W) \equiv \D(W,W) \,$ is 
identified with $\, \mbox{\rm End}_A(I) \,$.  On the other hand, 
if $\, W \,$ corresponds to the pair $\, (X , \mathcal{L}) \,$, then 
according to 
Proposition~\ref{dw}, $\,\D(W) \,$ is just the 
algebra $\, \D_{\mathcal L}(X) \,$ that interests us. 
In this way Theorem~\ref{chgr} translates any question about the 
algebras$\, \D_{\mathcal L}(X) \,$ into  a question about 
ideals in the Weyl algebra. It remains to give the translation
into these terms of the group $\, G \,$ of symmetries of $\, \gr \,$.
Note that if $\, \sigma \,$ is an automorphism of $\, A \,$ and 
$\, I \,$ ia  finitely generated torsion-free rank 1 $\, A $-module, 
then $\, \sigma_*(I) \,$ is a module of the same type: thus the 
automorphism group $\, \mbox{\rm Aut}(A) \,$ acts naturally on 
$\, \mathcal{R} \,$.
\begin{theorem}
\label{aut}
Under the bijection $\, \alpha \,$, the action of the group 
$\, \Gamma \,$ of KP flows corresponds to the action on 
$\, \mathcal{R} \,$ induced by the automorphisms 
$\, D \mapsto e^{p(z)} D e^{-p(z)} \,$ of $\, A \,$; while the map  
$\, \varphi \,$ corresponds to the map on $\, \mathcal{R} \,$ 
induced by the formal Fourier transform  
$\, (z \mapsto \partial, \, \partial \mapsto -z) \,$ of $\, A \,$.
\end{theorem}

Now, if $\, \sigma \,$ is an automorphism of (any 
algebra) $\, A \,$, and $\, M \,$ is any $\, A $-module, then it is 
trivial that $\, \mbox{\rm End}_A(M) \simeq \mbox{\rm End}_A(\sigma_{*}M) \,$.
Thus Theorem~\ref{aut} makes the ``if'' part of  Theorem~\ref{mainth}(i)
transparent.
\subsection*{Notes}
1. According to Dixmier (see \cite{D}), the automorphisms mentioned
in Theorem~\ref{aut} generate the full automorphism group of 
$\, A \,$; thus we may identify our symmetry group $\, G \,$ with 
$\, \mbox{\rm Aut}(A) \,$.
\vspace{0.1cm} \\
2. There are two routes available to prove the ``only if'' part of 
Theorem~\ref{mainth}(i).  If we use Dixmier's theorem, we can 
simply note that it translates into a  
known theorem of Stafford (see \cite{St}): if  $\, I \,$ and $\, J \,$ 
are two ideal classes of $\, A \,$, then their endomorphism rings are 
isomorphic (if and) only if $\, I \,$ and $\, J \,$  
belong to the same orbit of 
$\, \mbox{\rm Aut}(A) \,$ in $\, \mathcal{R} \,$.  Alternatively, after 
we have classified the orbits, this fact will follow from 
Theorem~\ref{lml} (whose proof in \cite{LM} does not use 
Stafford's theorem, nor Dixmier's).
\vspace{0.1cm} \\
3.  To get an idea of the depth of Stafford's theorem, let us give 
a proof (following \cite{CH3}) of a crucial special case: if $\, I \,$ 
is an ideal of $\, A \,$ whose endomorphism ring is isomorphic to  
$\, \mbox{\rm End}_A(A) = A\,$, 
then $\, I \simeq A \,$.  Let $\, (X, \mathcal L) \,$ be the pair 
corresponding to $\, I \,$; then $\, \D_{\mathcal L}(X) \,$ 
is isomorphic to $\, A \,$, hence $\, \O(X) \,$ 
is isomorphic to 
a mad subalgebra of $\, A \,$.  Another (nontrivial) theorem of 
Dixmier (see \cite{D}) says that all the mad subalgebras of $\, A \,$ 
are isomorphic to $\, \c[z] \,$; hence $\, X \simeq \a^1 \,$ 
and $\, \mathcal{L} \,$ is the trivial line bundle (because this 
is the only rank 1 torsion-free sheaf over $\, \a^1 $).  
According to Theorem~\ref{chgr}, it follows that $\, I \simeq A \,$.
The general case of Stafford's theorem is a relatively formal 
consequence of this special case (see \cite{St}, Corollary 3.2).
\vspace{0.1cm} \\
4. If we introduce the category $\, \mathfrak{P} \, $ with objects 
the primary decomposable subspaces of $\, \c[z] \,$ and morphisms 
$\, \D(V,W) \, $, then we could summarize Theorem~\ref{ch} 
by saying that we have an {\it equivalence of categories} between 
$\, \mathfrak{P} \, $ and the category of ideals in $\, A \,$ 
(regarded as a full subcategory of the category of right $\, A $-modules).
\vspace{0.1cm} \\
5.  Theorems~\ref{ch} and \ref{chgr} remain true ({\it mutatis mutandis}) 
if we replace the Weyl algebra by the ring of differential operators on any 
smooth affine curve (see \cite{CH1}).  Using this fact, we can 
sketch a proof of Theorem~\ref{smoothm}. Suppose that 
$\, X \,$ and $\, Y \,$ are smooth affine curves such that 
$\, \D(X) \,$ is Morita equivalent to $\, \D(Y) \,$.
Since these are Noetherian domains, that means that $\, \D(Y) \,$ 
is isomorphic to the endomorphism ring of some ideal in 
$\, \D(X) \,$, and hence to   $\, \D(V) \,$ for some 
primary decomposable subspace $\, V \,$ 
of $\, \O(X) \,$.  This in turn is 
isomorphic to some ring $\, \D_{\mathcal L}(X') \,$, where 
$\, X' \,$ is a curve with bijective normalization $\, X \to X' \,$. 
If $\, Y \,$ is not isomorphic to $\, \a^1 \,$, then
Theorem~\ref{mak} shows that $\, \D(Y) \,$ has only one mad 
subalgebra.  The same is therefore true of  
$\, \D_{\mathcal L}(X') \,$; extracting these mad subalgebras 
gives $\, \O(Y) \simeq \O(X') \,$, hence 
$\, Y \simeq X' \,$.  Since $\, Y \,$ is smooth, this implies $\, X = X' \,$, 
hence $\, Y \simeq X \,$. 
Finally, if $\, Y \,$ is isomorphic to $\, \a^1 \,$, then 
$\, \D(Y) \,$, and hence also  
$\, \D_{\mathcal L}(X') \,$, 
has more than one mad subalgebra, so Lemma~\ref{bisp} implies that 
$\, X \simeq \a^1 \,$.
\vspace{0.1cm} \\
6. Theorem~\ref{ch} is proved in \cite{CH1}; Theorem~\ref{aut} is proved 
in \cite{BW2}. 
\vspace{0.1cm} \\
7. A different view of the construction of Cannings and Holland,
and some further generalizations, can be found in \cite{BGK2}.

\section{The Calogero-Moser spaces}
\label{cm}
Our third realization of $\, \gr \,$ involves the 
{\it Calogero-Moser spaces} $\, \mathcal{C}_n \,$.  For each 
$\, n \geq 0 \,$, let $\, \tilde{\mathcal{C}}_n \,$ be the space of 
pairs $\, (X,Y) \,$ of complex $\, n \times n \,$ matrices such that 
$$
[X,Y] + I \ \text{has rank} \ 1 \, ,
$$
and let $\, \mathcal{C}_n := {\tilde{\mathcal{C}}_n}/{\mbox{\rm GL}(n,\c)} \,$, 
where the action of $\, g \in \mbox{\rm GL}(n,\c) \,$ is by simultaneous 
conjugation: $\, (X,Y) \mapsto (gXg^{-1}, gYg^{-1}) \,$.  One can show that 
$\, \mathcal{C}_n \,$ is an smooth irreducible affine variety of dimension 
$\, 2n \,$ ($\, \mathcal{C}_0 \,$ is supposed to be a point).
\begin{theorem}
\label{beta}
There is a natural bijection
$$
\beta \ : \ \mathcal{C} \, := \, \bigsqcup_{n \geq 0} \mathcal{C}_n \to \gr
$$
such that \\
(i) the action of $\, \Gamma \,$ on $\, \gr \,$ corresponds to the maps 
$\, (X,Y) \mapsto (X + p'(Y),Y) \,$ on $\, \mathcal{C}_n \,$; \\
(ii) the action of $\, \varphi \,$ on $\, \gr \,$ corresponds to the map 
$\, (X,Y) \mapsto (-Y,X) \,$ on $\, \mathcal{C}_n \,$; \\
(iii)  the action of the group $\, G \,$ on each $\, \mathcal{C}_n \,$ 
is transitive.
\end{theorem}

It follows from part (iii) of this Theorem that the spaces 
$\, \beta(\mathcal{C}_n) \,$ 
are the orbits of $\, G \,$ in $\, \gr \,$. To complete the proof 
of Theorem~\ref{mainth} we have only to check that $\, \beta^{-1} (W_n) \,$ 
belongs to $\, \mathcal{C}_n \,$: that is done in Example~\ref{dn} below.

The decomposition of $\, \gr \,$ in Theorem~\ref{beta} was originally 
obtained using ideas from the theory of integrable systems (see \cite{W2}).
Here we sketch a different method.  In view of Theorem~\ref{chgr}, it is 
enough to see why the space $\, \mathcal{R} \,$ of ideals in the 
Weyl algebra should decompose into the finite-dimensional 
spaces $\, \mathcal{C}_n \,$.  That can be understood by analogy with the 
corresponding commutative problem, namely, to describe the space 
$\, \mathcal{R}_0 \,$ of isomorphism classes of ideals in 
$\, A_0 := \c[x,y] \,$.  This problem is easy, because 
each ideal class in $\, A_0 \,$ has a unique 
representative of finite codimension; hence $\, \mathcal{R}_0 \,$ 
decomposes into the disjoint union of the {\it point Hilbert schemes} 
$\, \mbox{\rm Hilb}_n(\a^2) \,$ (that is, the spaces of ideals 
of codimension $\, n \,$) for $\, n \geq 0 \,$.  It is elementary that 
$\, \mbox{\rm Hilb}_n(\a^2) \,$ can be identified with 
the space of pairs $\, (X,Y) \,$ of commuting $\, n \times n \,$ 
matrices possessing a cyclic vector (see \cite{N}, 1.2); thus 
$\, \mbox{\rm Hilb}_n(\a^2) \,$ is the commutative analogue 
of the Calogero-Moser space $\, \mathcal{C}_n \,$. Because the Weyl 
algebra has no nontrivial ideals of finite codimension, it is not 
immediately clear how to adapt this discussion to the noncommutative 
case; however, there is a less elementary point of view which 
generalizes more easily. We may regard an ideal of $\, A_0 \,$ as a 
rank 1 torsion-free sheaf over $\, \a^2 \,$; it has a unique 
extension to a torsion-free sheaf over the projective plane 
$\, \P^2 $ trivial over the line at infinity. 
The classification of ideals by pairs of matrices 
can then be regarded as  (trivial) special case of Barth's  classification 
of framed bundles (of any rank) over $\, \P^2 \,$ 
(see \cite{N}, Ch.\ 2).  In a similar way, an ideal of the Weyl 
algebra determines a rank 1 torsion-free sheaf over a suitably defined 
quantum projective plane   $\, \P^2_q \,$; these can then be 
classified much as in the commutative case.
\subsection*{Notes}
1. Let us try to give something of the flavour of the noncommutative 
projective geometry needed to carry out the plan sketched above
(see, for example \cite{A}, \cite{AZ} for more details).  Let 
$\, X \subseteq \P^N \, $ be a projective variety, and let 
$\, \boldsymbol{A} = \oplus_{k \geq 0} A_k \,$ be its (graded) 
homogeneous coordinate ring. To any quasicoherent sheaf 
$\, \mathcal{M} \,$ over $\, X \,$ we can assign the graded 
$\, \boldsymbol{A}$-module 
$$
\boldsymbol{M} := \bigoplus_{k \in \Z} H^0(X, \mathcal{M}(k)) \, .
$$
A theorem of Serre (see \cite{S}) states that this defines an equivalence 
between the category of quasicoherent sheaves over $\, X \,$ and a certain 
quotient of the category of graded $\, \boldsymbol{A}$-modules 
(we have to divide out by the so-called {\it torsion modules}, in which 
each element is killed by some $\, A_k \,$).  
Thus many results about projective varieties can be formulated in 
a purely algebraic way, in terms of graded $\, \boldsymbol{A}$-modules; 
in this form the theory makes sense also for a noncommutative graded 
ring  $\, \boldsymbol{A} \,$.  The coordinate ring of the space 
$\, \P^2_q \,$ referred to above is the ring of noncommutative 
polynomials in three variables $\, x,y,z \,$ of degree 1\,, where 
$\, z \,$ commutes with everything, but $\, [x,y] = z^2 \,$.  It turns 
out that the homological properties of this ring are similar to those 
of the commutative graded ring $\, \c[x,y,z] \,$; in particular, the 
classification of bundles (of any rank) over $\, \P^2_q \,$ 
is similar to that of bundles over $\, \P^2 $ (see \cite{KKO}).
\vspace{0.1cm} \\
2.  The idea of using $\, \P^2_q \,$ to classify the ideals in 
the Weyl algebra is due to L.\ Le Bruyn (see \cite{LeB}).  However, Le 
Bruyn's chosen extension of an ideal in $\, A \,$ to a sheaf over 
$\, \P^2_q \,$ was in general not trivial over the line at 
infinity, so he did not obtain the decomposition of $\, \mathcal{R} \,$ 
into the Calogero-Moser spaces.  That was done in \cite{BW3} and 
(in a different way) in \cite{BGK1}.
\vspace{0.1cm} \\
3.  The connection between the spaces 
$\, \mbox{\rm Hilb}_n(\a^2) \,$ and 
$\, \mathcal{C}_n \,$ is actually much closer than we have indicated: 
$\, \mbox{\rm Hilb}_n(\a^2) \,$ is a hyperk\"ahler variety, 
and $\, \mathcal{C}_n \,$ is obtained by deforming the complex 
structure of $\, \mbox{\rm Hilb}_n(\a^2) \,$ within the 
hyperk\"ahler family.  See \cite{N}, Ch.\ 3, especially 3.45.
\vspace{0.1cm} \\
4. The assertions (i) and (ii) in Theorem~\ref{beta} are proved in 
\cite{BW2} (using the original construction of $\, \beta $), and 
in \cite{BW3} (using the construction sketched above). The fact that 
the two constructions agree is also proved in \cite{BW3}. 
\vspace{0.1cm} \\
5. Parts (i) and (ii) of Theorem~\ref{beta} reduce the proof 
of part (iii) (transitivity of the $\, G$-action) to an exercise in 
linear algebra. Unfortunately, the exercise seems to be quite difficult,
and the published solution in \cite{BW2} strays 
outside elementary linear algebra at one point (see 
Lemma 10.3 in \cite{BW2}).  P. Etingof has kindly pointed out to us 
that transitivity also follows easily from the fact that 
the functions $\, (X,Y) \mapsto \mbox{\rm tr} (X^k) \,$ and 
$\, (X,Y) \mapsto \mbox{\rm tr} (Y^k) \,$ generate 
$\, \O(\mathcal{C}_n) \,$ as a Poisson algebra 
(see \cite{EG}, 11.33).
\vspace{0.1cm} \\
6. In \cite{BW3}, Section 5 we have given an elementary construction 
of the map $\, \mathcal{R} \to \mathcal{C} \,$, in a similar spirit to the 
elementary treatment of the commutative case. It turns out that the inverse 
map  $\, \mathcal{C} \to \mathcal{R} \,$ can also be written down 
explicitly, as follows.  Let $\, (X,Y) \in \mathcal{C}_n \,$, and choose 
column and row vectors $\, v, w \,$ such that $\, [X,Y] + I = vw \,$. 
Define\footnote{We get this formula by combining Remark 5.4 in \cite{BW3}
with formula (3.5) in \cite{W2}.}
$$
\kappa \,:=\, 1 - w(Y - zI)^{-1} (X - \partial I)^{-1}v
$$
(thus $\, \kappa \,$ belongs to the quotient field of the Weyl 
algebra $\, A \,$).  Then the (fractional) right ideal 
\begin{equation}
\label{rep}
\det(Y - zI) \, A \,+\, \kappa \det(X - \partial I) \, A 
\,\subset\, \c(z)[\partial]
\end{equation}
represents the class in $\, \mathcal{R} \,$ corresponding to 
$\, (X,Y) \,$.  Using these formulae, it is possible to give a completely 
elementary proof that $\, \mathcal{R} \,$ decomposes into the 
spaces $\, \mathcal{C}_n \,$.  More details will appear elsewhere.

\section{The invariant $\, n \,$}
\label{n}

Theorem~\ref{beta} assigns to each $\, W \in \gr \,$ a 
non-negative integer $\, n \,$, namely, the index of the ``stratum'' 
$\, \mathcal{C}_n \,$ containing $\, \beta^{-1} (W) \,$.  Using 
Proposition~\ref{wxl} and Theorem~\ref{chgr}, 
we may equally well regard $\, n \,$ as an invariant of 
a pair $\, (X, \mathcal L) \,$, or of an ideal (class) in the Weyl 
algebra  $\, A \,$.  In this section we discuss various  
descriptions of this invariant.  The first two begin with an ideal class 
in $\, A \,$.

\subsection*{$ n \,$ as a Chern class}

We return to the quantum projective plane  $\, \P^2_q \,$ 
explained at the end of Section~\ref{cm}.  Let $\, M \,$ be an ideal 
class of $\, A \,$, and let $\, \mathcal{M} \,$ denote its 
unique extension to a sheaf over $\, \P^2_q \,$ 
trivial over the line at infinity.
Then we claim that
\begin{equation}
\label{h1}
n = \dim_{\c} H^1( \P^2_q , \mathcal{M}(-1)) \, .
\end{equation}
To see that, we need to give more details of the construction of the map 
$\, \mathcal{R} \to  \mathcal{C} \,$.  Recall that the homogeneous 
coordinate ring of $\, \P^2_q \,$ has three generators 
$\, x,y,z \,$. It turns out that multiplication by $\, z \,$ induces 
an isomorphism
$$
 H^1( \P^2_q , \mathcal{M}(-2)) \to 
H^1( \P^2_q , \mathcal{M}(-1) := V \, .
$$
If we use this isomorphism to identify these spaces, then 
multiplication by $\, x \,$ and $\, y \,$ gives us a pair $\, (X,Y) \,$ 
of endomorphisms of $\, V \,$: this is the point of $\, \mathcal{C} \,$ 
associated with $\, M \,$.  Obviously, the size of the matrices 
$\, (X,Y) \,$ is given by (\ref{h1}).
\subsection*{Note}
By analogy with the commutative case (see \cite{N}, Ch.\ 2), 
we would like to interpret 
$\, n \,$ as the second Chern class $\, c_2(\mathcal{M}) \,$.  However, 
at the time of writing, Chern classes have not yet been discussed in 
noncommutative projective geometry.
\subsection*{$ n \,$ as a codimension}

Again, let $\, M \,$ be an ideal of $\, A \,$.  By \cite{St}, Lemma 4.2, 
we may suppose that $\, M \,$ intersects $\, \c[z] \subset A \,$ 
nontrivially; let $\, I \,$ be the ideal in $\, \c[z] \,$ 
generated by the leading coefficients of the operators in 
$\, M \,$, and let $\, p(z) \,$ be a generator of $\, I \,$. Then 
$\, p^{-1}M \subset \c(z)[\partial] \,$ is a fractional ideal 
representing the class of $\, M \,$. 
Define a map $\, D \mapsto D_{+} \,$ from $\, \c(z)[\partial] \,$ 
to $\, A \, $ by 
$$
\Bigl( \sum_i f_i \partial^i \Bigr)_{+} = \, \sum_i (f_i)_{+} \partial^i \, ;
$$
here $\, f_{+} \,$ denotes the polynomial part of a rational 
function $\, f \,$ (that is, the polynomial such that 
$\, f - f_{+} \,$ vanishes at infinity).
Then we claim that
$$
n \ \, \text{\it is the codimension of} \,  \ 
(p^{-1}M)_{+} \, \ \text{\it in} \, \ A \, .
$$
A proof can be found in \cite{BW3}, Section 6, 
where it is shown that the quotient 
space $\, A/{(p^{-1}M)_{+}}$ can be identified with the (\v{C}ech) 
cohomology group on the right of (\ref{h1}).
\subsection*{Note}
The special representative for an ideal class that we used in this 
subsection is the same one as is given by the formula (\ref{rep}). 
It is the unique representative of the form 
$\, \D(\c[z], W) \,$ with $\, W \in \gr \,$ (cf.\ Theorem~\ref{chgr}).
\subsection*{The differential genus of a framed curve}

The following characterization of $\, n \,$ was one of 
the main results of \cite{W2}.
\begin{theorem}
\label{celldim}
Let $\, W \in \gr \,$.  Then the integer $\, n \,$ that we have 
associated to $\, W \,$ is equal to the 
dimension of the open cell in $\, \gr $ containing  $\, W \,$.
\end{theorem}

This theorem leads easily to a simple formula for 
calculating $\, n \,$ in concrete examples (cf.\ \cite{PS}, 7.4). 
Recall from Section~\ref{grad} that $\, W \,$ is constructed from a family 
of $\, \lambda$-primary subspaces $\, V_{\lambda} \subseteq \c[z] \,$ 
(one for each $\, \lambda \in \c \,$, and almost all of them equal 
to $\, \c[z] \,$).  In terms of these $\, V_{\lambda} \,$, we can 
calculate $\, n \,$ as follows.  First, we have 
$\, n = \sum_{\lambda} n_{\lambda} \,$, where $\, n_{\lambda} \,$ 
depends only on $\, V_{\lambda} \,$ (and is zero if 
$\, V_{\lambda} = \c[z] \,$).  To find $\, n_{\lambda} \,$, let
\begin{equation}
\label{1}
r_0 < r_1 < r_2 < \ldots
\end{equation}
be the numbers $\, r \,$ such that $\, V_{\lambda} \,$
contains a polynomial that vanishes exactly to order 
$\, r \,$ at $\, \lambda \,$.  For large $\, i \,$
we have $\, r_{i} = g + i \,$, where $\, g \,$ is the number of 
``gaps'' (non-negative integers that do not occur) in the sequence 
(\ref{1}).  Then we have
\begin{equation}
\label{2}
n_{\lambda} \  = \ \sum\limits^{}_{i \geq 0} \, (g + i  - r_i)\ .
\end{equation}
\begin{example}
\label{dn}
For the $\, 0$-primary space $\, V := \O(X_n) \,$ 
defined by (\ref{Xn}), 
the sequence (\ref{1}) is
$$
0 < n+1 < n+2 < \ldots \, ,
$$
whence $\, g = n \,$, and the right hand side of (\ref{2}) 
is equal to $\, n \,$.  
\end{example}

This calculation completes the proof of Theorem~\ref{mainth}(ii), and 
shows that we can identify the number $\, n \,$ associated with 
a pair $\, (X, \mathcal L) \,$ with the differential genus 
$\, d_{\mathcal{L}}(X) \,$ introduced in Section~\ref{grad}.
\begin{example}
If $\, Y_r \,$ is the curve with coordinate ring 
$\, \O(Y_r) := \c[z^2, z^{2r+1}] \,$ , 
then again $\, \O(Y_r) \,$ is $\, 0$-primary, and 
the sequence (\ref{1}) is
$$
0 < 2 < 4 <  \ldots < 2r < 2r+1 < \ldots \, .
$$
Hence $\, g = r \,$, and $\, d(Y_r) = 
r + (r-1)+ \ldots + 2 + 1 = r(r+1)/2 \,$.
\end{example}

In particular, $\, d(Y_2) = 3 \,$ so 
$\, Y_2 \,$ is differentially isomorphic to $\, X_3 \,$,
in agreement with G. Letzter (see \cite{L}).
\begin{example}
\label{34}
Here is the simplest example to show that in general 
$\, d_{\mathcal{L}}(X) \,$ depends on $\, \mathcal{L} \,$, not 
just on $\, X \,$.
Let $\, V \,$ be the $\, 0$-primary space spanned by 
$\, \{ z^i : i \not= 2,3 \} \,$.  Then the the sequence (\ref{1}) is
$$
0 < 1 < 4 < 5 < \ldots \, ,
$$
whence $\, n = 4 \,$.  Clearly, $\, V \,$ is a maximal module 
over the ring $\, \O(X_3) \,$, and thus corresponds to 
a maximal torsion-free (but not locally free) sheaf $\, \mathcal{L} \,$ 
over $\, X_3 \,$.  For this sheaf $\, \mathcal{L} \,$ 
we therefore have $\, d_{\mathcal{L}}(X_3) = 4\,$, 
and the ring $\, \D_{\mathcal{L}}(X_3) \,$ is isomorphic 
to $\, \D(X_4) \,$.
\end{example}

\subsection*{The Letzter-Makar-Limanov invariant}

Next, we describe the invariant originally used in \cite{LM} to 
distinguish the rings $\, \D(X_n) \,$.  We 
return temporarily to the case of any affine curve $\, X \,$, 
with normalization $\, \tilde X \,$ and function field $\, \k \,$; 
as usual (see Proposition~\ref{conc}), 
we view $\, \D(X) \,$ and $\, \D(\tilde X) \,$ 
as subalgebras of 
$\, \D(\k) \,$.  In general, $\, \D(X) \,$ is not 
contained in $\, \D(\tilde X) \,$; however, the associated 
graded algebra $\, \makebox{\rm gr} \, \D(X) \,$ is always 
contained in $\, \makebox{\rm gr} \, \D(\tilde X) \,$ 
(see \cite{SS}, 3.11).  
In the case that most concerns us when $\, \tilde X = \a^1 \,$, 
this simply means that the {\it leading} coefficient of each operator 
in $\, \D(X) \,$ is a polynomial (although the other 
coefficients may be rational functions, as we saw in Example \ref{x1}). 
Continuing Theorem \ref{six}, we have 
\begin{theorem}
\label{cdfinite}
Each of the conditions in Theorem \ref{six} is equivalent to: 
$$
\makebox{\rm gr} \, \D(X) \, \ 
\text{has finite codimension in} \ 
\, \makebox{\rm gr} \, \D(\tilde X) \, .
$$
\end{theorem}

In our case, when $\, \tilde X = \a^1 \,$ 
and $\, X \,$ is a framed curve,
$\, \makebox{\rm gr} \, \D(X) \,$ is a subalgebra of 
finite codimension in $\, \c[z, \zeta ] \,$; we call its codimension 
the {\it Letzter-Makar-Limanov invariant} of $\, X \,$, and denote 
it by $\, LM(X) \,$.  The definition of $\, LM(X) \,$ uses the 
standard filtration on $\, \D(X) \,$; nevertheless, in 
\cite{LM} it is proved that it depends only on the isomorphism class 
of the algebra $\, \D(X) \,$; that is, if $\, X \,$ and 
$\, Y \,$ are differentially isomorphic framed curves, then 
$\, LM(X) =  LM(Y) \,$.  On the other hand, it is not hard to 
calculate that $\, LM(X_n) = 2n \,$ (see \cite{LM}, Section 5).  
Combined with Theorem~\ref{mainth}, that gives
\begin{theorem}
Let $\, X \,$ be any framed curve. Then  $\, 2 \, d(X) = LM(X) \,$.
\end{theorem}
\subsection*{Notes}
1.  Theorem ~\ref{cdfinite} is proved (though 
not explicitly stated) in \cite{SS}, 3.12.
\vspace{0.1cm} \\
2.  In \cite{LM} the rings $\, \D_{\mathcal{L}}(X) \,$ 
(for $\, \mathcal{L} \not= \O_{X} \,$) are not considered; 
however, it is not hard to extend the discussion to include that 
case.  Thus we can define the invariant $\, LM(\D(W)) \,$ for any 
$\, W \in \gr \,$, and Theorem~\ref{mainth} shows that it is 
equal to $\, 2n \,$.
\vspace{0.1cm} \\
3. It is possible to prove directly 
(that is, without using Theorem~\ref{mainth}) that $\, LM(\D) \,$ 
is twice the number $\, n \,$ defined by (\ref{h1}). The interested 
reader may see \cite{B}.
\vspace{0.3cm}

All our descriptions of $\, n \,$ so far have been specific to our 
particular situation.  It is natural to ask whether $\, n \,$ is a 
special case of some general invariant of rings that is able 
to distinguish between different Morita equivalent domains.  Our last 
two subsections are attempts in that direction.
\subsection*{Pic and Aut}

Let $\, \D \,$ momentarily be any domain (associative  
algebra without zero divisors) over $\, \c \,$. The following 
idea for obtaining subtle invariants of the isomorphism class of 
$\, \D \,$ is due to Stafford (see \cite{St}).  
Consider the group\footnote{More properly, 
we should write $\, \mbox{\rm Pic}_{\c}(\D) \,$ to indicate that we 
consider only equivalences that commute with multiplication by scalars.
For a similar reason, we should write $\, \mbox{\rm Aut}_{\c} \,$ too.} 
$\, \mbox{\rm Pic}(\D) \,$ 
of all Morita equivalences of $\, \D \,$ with itself, that is, of all 
self-equivalences of the category Mod-$ \D \,$ of 
(say right) $\, \D$-modules.  
Each such equivalence is given by tensoring with a suitable 
$\, \D$-bimodule, so we may also think of 
$\, \mbox{\rm Pic}(\D) \,$ 
as the group of all invertible $\,\D $-bimodules.  Each automorphism 
of $\, \D \,$ induces a self-equivalence of  
Mod-$ \D \,$, so there is a natural map
\begin{equation}
\label{omega}
\omega : \mbox{\rm Aut}(\D) \to \mbox{\rm Pic}(\D) \, .
\end{equation}
Although the group $\, \mbox{\rm Pic}(\D) \,$ is a Morita invariant 
of $\, \D \,$, the automorphism group and the map 
$\, \omega \,$ are not.  

We return to our case, where $\, \D \,$ is one of the 
algebras $\, \mbox{\rm End}_A(I) \,$ (or $\, \D_{\mathcal{L}}(X) \,$).  
In general, the kernel of $\, \omega \,$ 
consists of the {\it inner} automorphisms of $\, \D \,$; in 
our case these are trivial, so $\, \omega \,$ is injective.  For the 
Weyl algebra $\, A \,$, Stafford showed that $\, \omega \,$ is an isomorphism. 
We thus have a natural inclusion
$$
\mbox{\rm Aut}(\D) \hookrightarrow 
\mbox{\rm Pic}(\D) \simeq \mbox{\rm Pic}(A) = 
\mbox{\rm Aut}(A) 
$$
(the isomorphism from $\, \mbox{\rm Pic}(\D) \,$ to 
$\, \mbox{\rm Pic}(A) \,$ is defined by tensoring with the 
$\, \D$-$A$-bimodule $\, I \,$).  
Recalling that the group $\, \mbox{\rm Aut}(A) \,$ 
acts transitively on $\, \mathcal{C}_n \,$, one can calculate that 
the isotropy group of the point in $\, \mathcal{C}_n \,$ corresponding to 
$\, I \,$ is exactly this subgroup 
$\, \mbox{\rm Aut}(\D) \,$.  It follows that we have 
a natural bijection
$$
\mathcal{C}_n \simeq 
{\mbox{\rm Pic}(\D)}/{\mbox{\rm Aut}(\D)}\ ,
$$
so it is tempting to claim that our invariant $\, n \,$ is given by
\begin{equation}
\label{quotient}
2n = \dim_{\c}\,{\mbox{\rm Pic}(\D)}/{\mbox{\rm Aut}(\D)}\ .
\end{equation}
The flaw in this is that the structure of algebraic variety on the 
quotient ``space'' in (\ref{quotient}) has been imposed {\it a posteriori}, 
and has not been extracted intrinsically from the algebra 
$\, \D \,$.
\subsection*{Note}
In view of the above, we may hope that there should be 
(at least for some algebras $\, \D \,$) a natural structure 
of (infinite-dimensional) algebraic group on 
$\, \mbox{\rm Pic}(\D) \,$ for which 
$\, \mbox{\rm Aut}(\D) \,$ would be a closed subgroup.  In our 
case, we can identify $\, \mbox{\rm Pic}(\D) \,$ with 
$\, \mbox{\rm Aut}(A) \,$, which does indeed have a natural 
structure of algebraic group; however, for this structure 
$\, \mbox{\rm Aut}(\D) \,$ is not a closed subgroup (see 
\cite{BW2}, Section 11 for more details).
\subsection*{Mad subalgebras}

The idea behind our final description of $\, n \,$ is very simple, 
namely: $\, n \,$ should measure the ``number'' of mad subalgebras of  
$\, \D(X) \,$.  Let us formulate a precise statement.
For each $\, W \in \gr \,$ with invariant $\, n \,$,
 we may choose an isomorphism 
$$
\phi : \D(W)  \to \D(X_n) \, .
$$
Since $\, \D^{0}(W) \,$ is a mad subalgebra of 
$\, \D(W) \,$,  $\, B := \phi(\D^{0}(W)) \,$ 
is a mad subalgebra of $\, \D(X_n) \,$. 
Furthermore, $\, \phi \,$ extends to an isomorphism of quotient 
fields, in particular, it maps $\, z \in \c(z)[\partial] \,$ 
to some element $\, u := \phi(z) \,$ in the quotient field of $\, B \,$.  
Clearly, $\, \c[u] \,$ is the integral closure of $\, B \,$.  
According to \cite{LM}, the integral closure $\, \overline{B} \,$ of {\it any} 
mad subalgebra $\, B \,$ is isomorphic to $\, \c[u] \,$: 
we shall call a choice 
of generator for  $\, \overline{B} \,$  a {\it framing} of $\, B \,$. 
Thus the above isomorphism $\, \phi \,$ gives us a framed mad subalgebra 
$\, (B, u) \,$ of $\, \D(X_n) \,$. Any two choices of 
$\, \phi \,$ differ only by an automorphism of $\, \D(X_n) \,$,
so the {\it class} (modulo the action of 
$\, \mbox{\rm Aut} \,\D(X_n) $) of the 
framed mad subalgebra we have obtained 
depends only on $\, W \,$.  Moreover (cf.\  Section~\ref{grad}, Note 1),
if we replace $\, W \,$ by 
$\, gW \,$, where $\, g \,$ belongs to the group $\, \Gamma \,$ of 
KP flows, then conjugation by $\, g \,$ defines an
isomorphism of $\, \D(gW) \,$ with $\, \D(W) \,$ which is the identity 
on  $\,\D^{0} \,$, 
so the isomorphism
$$
\D(gW) \to \D(W) \stackrel{\phi}
{\longrightarrow} \D(X_n)
$$
defines the same framed mad subalgebra as $\, \phi \,$.  
It follows that we have constructed a well-defined map
\begin{equation}
\label{conjecture}
{\mathcal{C}_n}/{\Gamma} \to  
\{\text{classes of framed mad subalgebras in}\  \D(X_n) \} \, .
\end{equation}
\begin{theorem}
\label{latest}
The map (\ref{conjecture}) is a set-theoretical bijection.
\end{theorem}

We will explain the proof elsewhere.  
Since the (categorical) quotient $\, {\mathcal{C}_n}/\!/{\Gamma} \,$ 
is $\, n$-dimensional, we should like to interpret 
$\, n \,$ as the dimension of the ``space'' of (classes of framed) mad 
subalgebras of $\, \D(X_n) \,$.  However, the word  ``space'' 
here is open to even more serious objections than in the preceding subsection.
\subsection*{Notes}
1. In the definition of a framed mad subalgebra $\, (B, u) \,$ we did 
not assume {\it a priori} that the curve 
$\, \mbox{\rm Spec} \, B \,$ was free of 
multiple points (indeed, that was not proved in \cite{LM}).  This 
momentary inconsistency of terminology is resolved by Theorem~\ref{latest}, 
which asserts ({\it inter alia}) that every mad subalgebra $\, B \,$ 
arises from the construction described above; in particular, that 
$\, \mbox{\rm Spec} \, B \,$ is a framed curve as defined earlier.
\vspace{0.1cm} \\
2.  In the case $\, n = 0 \,$, the left hand side of 
(\ref{conjecture}) is a point, so Theorem~\ref{latest} becomes a well known
result of Dixmier: in the Weyl algebra there is only one class of 
mad subalgebras (see \cite{D}).  
\section{Higher dimensions}
\label{higher}

\subsection*{The examples of Levasseur, Smith and Stafford}

Let $\, \g \,$ be a simple complex Lie algebra, and let $\, O \,$ be 
the closure of the minimal nilpotent orbit in $\, \g \,$.  Let 
$\, \g = \mathfrak{n}_{-} \oplus \h \oplus \mathfrak{n}_{+}$ be 
a triangular decomposition of $\, \g \,$; then $\, O \cap \mathfrak{n}_{+} \,$
breaks up into several irreducible components $\, X_i \,$.  In 
\cite{LSS} it is shown that in some cases the ring $\, \D(X_i) \,$ can be 
identified with $\, U(\g)/J \,$, where $\, J \,$ is a certain distinguished 
completely prime primitive ideal of $\, U(\g) \,$ (the 
{\it Joseph ideal}).  The examples of differential isomorphism arise in 
the case $\, \g = \mathfrak{s}\mathfrak{o}(2n,\c) \,$ (with $\, n \geq 5 $), 
because in that case there are two nonisomorphic components $\, X_1 \,$ and 
$\, X_2 \,$ of this kind.  They can be described quite explicitly: 
$\, X_1 \,$ is the quadric cone $\, \sum z_i^2 = 0 \,$ in $\, \c^{2n-2} \,$, 
and $\, X_2 \,$ is the space of skew-symmetric $\, n \times n \,$ 
matrices of rank $\, \leq 2 \,$.  In contrast to what we saw for curves, 
these spaces $\, X_1 \,$ and $\, X_2 \,$ are quite different topologically.
\subsection*{Morita equivalence}
There are several papers that study differential 
{\it equivalence} in dimension $ > 1 \,$.  In view of Theorem~\ref{six}, 
attention has focused on the question of when 
a variety $\, X \,$ is differentially equivalent to its normalization 
$\, \tilde X \,$.  Of course, in dimension $ > 1 \,$ the normalization 
is not necessarily smooth:  in \cite{J1} there are 
examples of differential equivalence in which $\, \tilde{X}\,$ is not smooth 
(they can be thought of as  generalizations of the 
monomial curves of Example~\ref{wn}).
Another point that does not arise for curves is 
that the condition that $\, X \,$ be {\it Cohen-Macaulay} plays an 
important role (we recall that every curve is Cohen-Macaulay). For 
example, a theorem of Van den Bergh states that if 
$\, \D(X) \,$ is simple, then $\, X \,$ must be 
Cohen-Macaulay (see \cite{VdB}, Theorem 6.2.5).  For varieties with 
{\it smooth} normalization, there are good generalizations of 
at least some parts of Theorem~\ref{six}. For example, piecing together  
various results scattered through the literature, 
we can get the following.
\begin{theorem}
\label{HSCS}
Let $\, X \,$ be an (irreducible) affine variety with smooth 
normalization $ \tilde{X}\,$.  Then the following are equivalent.
\begin{enumerate}
\item The normalization map $\, \pi : \tilde X \to X \,$ is bijective and 
$ X $ is Cohen-Macaulay.
\item The algebras $\, \D(\tilde X) \,$ and $\, \D(X) \,$ 
are Morita equivalent.
\item The ring $\, \D(X) \,$ is simple.
\end{enumerate}
\end{theorem}

Beautiful examples are provided by the varieties of 
{\it quasi-invariants} of finite reflection groups 
(see \cite{BEG}, \cite{BC}): 
here $ \tilde{X}\,$ is the affine space $\, \a^m $, so these 
examples are perhaps the natural higher-dimensional generalizations 
of our framed curves.
\subsubsection*{References for the proof of Theorem~\ref{HSCS}}
For the implication ``(1) $ \Rightarrow $ (2)'' in Theorem~\ref{HSCS} we 
are relying on the recent preprint \cite{BN} (at least 
in dimension $\, > 2 \,$: 
for surfaces it was proved earlier in \cite{HS}).  For the rest, the 
implication ``(2) $ \Rightarrow $ (3)'' is trivial, and the fact that 
$\, \D(X) \,$ simple implies $\, X \,$ Cohen-Macaulay is 
the theorem of Van den Bergh mentioned above. The only remaining 
assertion in Theorem~\ref{HSCS} is that if $\, \D(X) \,$ is simple 
then $\, \pi \,$ is bijective.  Suppose $\, \D(X) \,$ is simple. 
Then by \cite{SS}, 3.3, $\, \D(X) \,$ 
is isomorphic to the endomorphism 
ring of the right $\, \D(\tilde X) $-module $ \, P \,$ defined 
by (\ref{P}); the dual basis lemma then implies that 
$\, P \,$ is a projective $\, \D(\tilde X) $-module. It now follows from  
\cite{CS}, Theorem 3.1 that $\, \D(X) \,$ is a 
{\it maximal order}, then 
from \cite{CS}, Corollary 3.4 that $\, \pi \,$ is bijective.
\subsection*{Non-affine varieties}
In this paper we have considered only affine varieties.  However, the problem 
of differential {\it equivalence} has an obvious generalization  
to arbitrary (for example, projective) varieties $\, X \,$. Namely: 
on $\, X \,$ we have the {\it sheaf} $\, \D_X \,$ of differential 
operators (whose sections over an affine open set 
$\, \mbox{Spec}\,A \,$ are the ring $\, \D(A) \,$), and given two 
varieties $\, X \,$ and $\, Y \,$, we can ask whether the categories 
of $\, \O $-quasicoherent sheaves of  modules over $\, \D_X \,$ and 
$\, \D_Y \,$ are equivalent.  For $\, X \,$ affine, the global section 
functor gives an equivalence between the categories of $\, \D_X$-modules 
and of $\, \D(X) $-modules, so we recover our original problem.
The available evidence (namely \cite{SS} and \cite{BN}) suggests that 
results about the affine case carry over to this more general situation.

The question of differential isomorphism does not make sense for 
sheaves; however, we can always consider the ring $\, \D(X) \,$ of 
global sections of $\, \D_X \,$, and ask when $\, \D(X) \,$ and 
$\, \D(Y) \,$ are isomorphic.  In general, $\, \D(X) \,$ may be 
disappointingly small: for example, if $\, X \,$ is a smooth projective 
curve of genus $\, > 1 \,$, then we have no global vector fields, so 
$\, \D(X) = \c \,$.  Probably the question is a sensible one 
only if  $\, X \,$ is close to being a 
$\, \D$-{\it affine} variety (for which the global section 
functor still gives an equivalence between  
$\, \D_X$-modules and $\, \D(X) $-modules). As far as we know, there are not 
yet any papers on this subject: however, the question of Morita 
eqivalence of rings  $\, \D(X) \,$ has been studied in \cite{HoS} 
(where  $\, X = \P^1 $, this being the only $\, \D$-affine smooth 
projective curve); and in \cite{J2} (where $\, X \,$ is a 
weighted projective space).  We should like to state one of the 
results of \cite{HoS}, since it is very close to our framed curves. 
Let $\, X \,$ be a ``framed projective curve'', that is, we have a 
bijective normalization map $\, \P^1 \to X \,$.  Then Holland and 
Stafford show that the rings $\, \D(X) \,$ (for $\, X \,$ singular) 
are all Morita equivalent to each other, but {\it not} to $\, \D(\P^1) \,$. 
A key point is that the although $\, \P^1 $ is $\, \D$-affine, the 
singular curves $\, X \,$ are not.
\medskip

\noindent
\small 
\textbf{Acknowledgments}.  The authors were partially supported by the 
National Science Foundation (NSF) grant DMS 00-71792 and an A.\ P.\ 
Sloan Research Fellowship.  The second author is grateful to the 
Mathematics Department of Cornell University for its hospitality 
during the preparation of this article.
\normalsize
\bibliographystyle{amsalpha}

\end{document}